\documentclass[12pt]{article}
\usepackage{amsmath,amsthm,amsfonts,amssymb,verbatim}
\newtheorem{thm}{Theorem}[section]
\newtheorem{lemma}[thm]{Lemma}
\newtheorem{defn}[thm]{Definition}
\newtheorem{cor}[thm]{Corollary}
\newtheorem{proposition}[thm]{Proposition}

\newtheorem{prop}[thm]{Proposition}

\theoremstyle{definition}

\newtheorem{remark}[thm]{Remark}
\theoremstyle{remark}
\DeclareMathOperator{\Fix}{Fix}
\DeclareMathOperator{\FC}{Fix_{c}}

\DeclareMathOperator{\Diff}{Diff}
\DeclareMathOperator{\MCG}{MCG}

\DeclareMathOperator{\Homeo}{Homeo}
\DeclareMathOperator{\Int}{int}
\DeclareMathOperator{\Per}{Per}

\newcommand{\cA}{{\cal A}}
\newcommand{\R}{\mathbb R}
\newcommand{\E}{{\cal E}}
\newcommand{\F}{{\cal F}}

\newcommand{\N}{{\cal N}}
\newcommand{\cP}{{\cal P}}
\newcommand{\cR}{{\cal R}}

\newcommand{\W}{{\cal W}}
\newcommand{\Z}{\mathbb Z}
\def\WF{\W}

\title{Fixed Points of abelian actions}
\author{John Franks,\thanks{Supported in part by NSF grant DMS0099640.}\ \ 
Michael Handel\thanks{Supported in part by NSF grant DMS0103435.}\ \
and Kamlesh Parwani \thanks{Supported in part by NSF grant DMS0244529.}}

\def\ti{\tilde}

\def\sl3z{SL(3, \mathbb Z)}

\def\G{{\cal G}}

\def\cS{{\cal S}}

\begin{document}
\maketitle
\begin{abstract}
We prove that if $\F$ is an abelian group of
$C^1$ diffeomorphisms isotopic to the identity
of a closed surface $S$ of genus at least two then
there is a common fixed point for all elements of $\F.$
\end{abstract}

\section{Introduction and Notation}

The Lefschetz theorem implies that if $S$ is a closed surface of
nonzero Euler characteristic and $f$ is an element of the group
$\Diff_0(S)$ of $C^1$ diffeomorphisms of $S$ isotopic to the identity,
then $f$ has a fixed point.  Several papers (see \cite{B1},\cite{B2},
\cite{han:commuting},\cite{FHP} \cite{F}) address the question of when
an abelian subgroup $\F$ of $\Homeo_0(S^2)$ or $\Diff_0(S^2)$ must
have a global fixed point, meaning a point that is fixed by every
element of the group.  The $S= S^2$ case is settled in \cite{FHP}
where it is shown that there is always a subgroup $\F_0$ of index at
most two with a global fixed point and that one can choose $\F_0 = \F$
if a certain $\Z_2$-invariant vanishes.  In this paper we therefore
restrict our attention to those surfaces $S$ with negative Euler
characteristic.

The fixed point set $\Fix(f)$ is partitioned into Nielsen classes,
each with its own Lefschetz index.  In our current context, there is a
preferred Nielsen class called the contractible Nielsen class whose
Lefschetz index is nonzero. (Full definitions and details are provided
in section~\ref{nielsen}.)  The set of contractible fixed points of
$f$ will be denoted $\FC(f)$ and if $\F$ is a subgroup of $\Diff_0(S)$
then $\FC(\F)$ will denote the points which are in $\FC(f)$ for all
$f$ in $\F.$ We can now state our main theorem.

\begin{thm}\label{thm:main}
Suppose $S$ is a closed oriented surface of genus at least two and
that $\F$ is an abelian subgroup of $\Diff_0(S)$
Then $\FC(\F)$ is non-empty.  In particular $\Fix(\F)$ is non-empty.
\end{thm}

There are two aspects of our proof that are of independent interest.  In section~\ref{rnc} we develop a theory of  Nielsen classes relative to certain infinite compact subsets $K$ of $S$ and more importantly of their Lefshetz indices.  The key   here is in isolating the property of $K$ that is required for the standard definitions to work.

In section~\ref{rc} we define reducing curves for a finite subset of $\F$ relative to certain $\F$-invariant compact sets.  These are analogous to the reducing curves that occur for Thurston canonical form, allowing us to decompose the surfce into simpler subsurfaces on which the group action is sufficiently well understood as to allow analysis.  The methods of \cite{FHP} are used in, but are not sufficient for, the construction of the set of reducing curves.

\section{Contractible Fixed Points} \label{nielsen}

Throughout this article $S$ will denote a closed oriented surface of
genus at least two.

A path $\alpha$ with endpoints in $\Fix(f)$ is an {\em
$f$-Nielsen path} if $f(\alpha)$ is homotopic to $\alpha$ rel
endpoints.  Elements of $\Fix(f)$ are {\em $f$-Nielsen equivalent} if
they are connected by an $f$-Nielsen path.  As the notation indicates,
this defines an equivalence relation on $\Fix(f)$ whose equivalence
classes are referred to as {\em $f$-Nielsen classes}.  When $f$ is
understood we refer to these simply as {\em Nielsen paths} and {\em
Nielsen classes}.

Since $S$ has arbitrarily small contractible neighborhoods, each
Nielsen class is an open subset of $\Fix(f)$, and since the complement
of a Nielsen class is a union of Nielsen classes, each Nielsen class
is also closed.  In particular, there are only finitely
many Nielsen classes.

If $\ti f : \ti S \to \ti S$ is a lift to the universal cover $\ti S$
of $S$ then any path $\ti \alpha$ with endpoints in $ \Fix(\ti f)$
projects to a Nielsen path in $S$. Conversely, if $\alpha$ is a
Nielsen path and $\ti f$ fixes one endpoint of a lift $\ti \alpha$ of
$\alpha$ then it also fixes the other endpoint of $\ti \alpha$.  Thus
$\Fix (\ti f)$ is either empty or projects onto a single Nielsen class
in $\Fix(f)$ that we say is {\em the Nielsen class determined by $\ti
f$}.  If $T$ is a covering translation of $\ti S$ then $\Fix(\ti f)$
and $\Fix(T \ti f T^{-1}) = T(\Fix(\ti f))$ project to the same
Nielsen class of $\Fix(f)$.  Conversely, if $\Fix(\ti f)$ and
$\Fix(\ti f')$ project to the same Nielsen class then there is a
covering translation $T$ such that $\ti f$ and $ T \ti f T^{-1}$ have
a common fixed point and hence are equal.

If $f$ is isotopic to the identity then there is a unique lift $\ti f
: \ti S \to \ti S$, called the {\em identity lift of $f$}, that
commutes with all covering translations of $\ti S$.  Equivalently, if
$f_t$ is any isotopy from $f_0 = $ identity to $f_1 = f$, and if $\ti
f_t$ is the lift of $f_t$ with $\ti f_0 = $ identity, then $\ti f_1$
is the identity lift.  The uniqueness of the identity lift follows
from the fact that the group of covering translations has trivial
center and so uses the fact that the Euler characteristic of $S$ is
negative.  The identity lift has non-empty fixed point set. Indeed, if
this fails then the unit tangent vector at $\ti x$ in the direction of the geodesic segment from $\ti x$ to $\ti f(\ti x)$
determines a non-zero equivariant vector field on $\ti S$ and hence a
non-zero vector field on $S$.  The Nielsen class determined by the
identity lift is called {\em the contractible Nielsen class} and is denoted $\FC(f)$. Elements of $\FC(f)$ are said to be {\em contractible}.  (The terminology comes from the fact that if $f_t$ is an isotopy
with $f_0 = id$ and $f_1 =f$ then  the closed loop $f_t(p)$ is contractible in $S.$) 
 Note that if $x \in \FC(f)$ and $\ti f$ is the identity
lift of $f$ then $\ti f$ fixes every lift of $x$ since it fixes some
lift of $x$ and commutes with every covering translation of $\ti S$.

Suppose now that $\F$ is a torsion free finitely generated abelian subgroup of $\Homeo_0(S)$.   Because the identity lift $\ti f$ of each $f \in \F$ is unique,  $f \to \ti f$ defines a lift  of $\F$ to a subgroup $\ti \F$ of $\Homeo(\ti S)$ called the {\em identity lift of $\F$}.  We say that $x \in S$ is a {\em global fixed point} if $x \in \Fix(f)$ for all $f \in \F$ and is a {\em global contractible fixed point} if $x \in \Fix_c(f)$ for all $f \in \F$; these sets are denoted $\Fix(\F)$ and $\Fix_c(\F)$ respectively.

\begin{lemma}  If $\{f_1,\ldots,f_m\}$ generate an abelian subgroup 
$\F$ of $\Homeo_0(S)$ then $\cap_{i=1}^m \FC(f_i) \subset \FC(\F)$.
\end{lemma}

\begin{proof}  If $\ti x \in \ti S$ is a lift of $x \in \cap_{i=1}^m \FC(f_i)$ then $\ti x \in \Fix(\ti f_i)$ for the  identity lift $\ti f_i$ of each $f_i$.  Since   the $\ti f_i$'s generate the identity lift $\ti \F$ of $\F$, we have $x \in \FC(\ti f)$ for each $\ti f \in \ti \F$. Thus   $x \in \FC(\F)$.
\end{proof}

Our interest in contractible fixed points comes from the following lemma.

\begin{lemma} \label{lem:contractible} 
If there exists   a global contractible fixed point  for a finite
index subgroup  of $\F$ then there exists   a global
contractible fixed point of $\F$.
\end{lemma}

\begin{proof}  If $\F_0$ has finite index in $\F$, if $x \in \FC(\F_0)$ and if $\ti \F_0$ is the identity lift of $\F_0$   then $\ti x \in \Fix(\ti F_0)$ for any lift $\ti x$ of $x$.   The orbit of $\ti x$ under $\ti \F$ is therefore finite and it follows from Theorem 1.2 of \cite{FHP} that $\Fix(\ti \F)$, and hence its projected image $\Fix_c(\F)$, are non-empty. 
\end{proof}

By combining Lemma~\ref{lem:contractible} with a result from \cite{FHP} we obtain a criterion for proving the existence of  global contractible fixed points.   

\begin{prop} \label{prop:lift} 
Suppose that $\F_0$ is a finite index subgroup of a finitely generated
abelian subgroup $\F$ of $\Homeo_0(S)$ and that $K$ is a compact
$\F_0$-invariant set.  Suppose also that $\gamma \subset S\setminus K$ is essential in $S\setminus K$ but inessential in $S$.   Finally, suppose  that there is a lift $\ti \gamma$ of $\gamma$ to the universal
covering $\ti S$  that   is preserved up to isotopy
rel $\ti K$, the complete pre-image of $K,$ by each $\ti f$ in the identity lift $\ti F_0$ of $F_0$.
  Then $\FC(\F) \ne \emptyset$.
\end{prop}
\begin{proof}
Let $\ti K_{\ti \gamma}$ be the intersection of $\ti K$ with the disk in $\ti S$ bounded
by $\ti \gamma.$  Then $\ti K_{\ti \gamma}$ is compact and is $\ti \F_0$-invariant, meaning that it is 
  $\ti f$-invariant  for every $\ti f \in \ti \F_0$.  
It follows from Theorem 1.2 of \cite{FHP} that  $\Fix(\ti \F_0) \ne \emptyset$ and hence that $\FC(\F_0) \ne \emptyset$.   
By Lemma~\ref{lem:contractible},  $\FC(\F)\ne \emptyset$.
\end{proof}

\section{Relative Nielsen Classes}  \label{rnc}

Some definitions in this section apply beyond the context of surface homeomorphisms so we consider   homeomorphisms $h : M \to M$ of a manifold of unspecified dimension.  For any compact $h$-invariant set $K$ we can define isotopy and homotopy rel   $K$ and in particular can define Nielsen classes rel $K$.  However, to make sense of the Lefschetz index of a Nielsen class we must place restrictions on $K$.    In particular we want the relative Nielsen classes to be open and closed subsets of $\Fix(h)$.  We make the following
definitions.

\begin{defn}
If $\alpha(s)$ and $\beta(s)$ are paths with the same
endpoints in a space $M$ and $K$
is a closed subset of $M$ we will say that $\alpha$ and $\beta$ are {\em homotopic rel $K$}
provided
there is a homotopy $\alpha_t(s),\  t \in [0,1]$
satisfying
\begin{enumerate}
\item  $\alpha_0 = \alpha$ and $\alpha_1 = \beta.$
\item  If  $\alpha_{t_0}(s) \in K$ for some $s$ and some $t_0$ 
then $\alpha_{t}(s) = \alpha_{t_0}(s)$ for all $t \in [0,1]$.
\item $\alpha_t(0) = x,\ \alpha_t(1) = y$ for all $0 \le t \le 1$.
\end{enumerate}

\end{defn}

Equivalently $\alpha_0(s)$ and $\alpha_1(s)$ are homotopic rel $K$
if they are homotopic as maps
$[0,1] \to S$ relative to the subset
$\alpha^{-1}(K) \cup \{0\} \cup \{1\}$ of $[0,1]$.

\begin{defn}
If $h: M \to M$ is a  homeomorphism with compact fixed point
set and $K$ is a compact invariant set we will say that an arc $\alpha$ 
joining points $x,y \in \Fix(h)$ is a
{\em Nielsen arc for $h$ rel $K$} if there is a homotopy rel $K$ from
$\alpha$ to $h(\alpha).$
\end{defn}

\begin{defn} 
If $x \in \Fix(h)$ and if there are arbitrarily small
neighborhoods $U_x \subset V_x$ of $x$ with the property that for
every $y \in U_x \cap \Fix(h)$ there is a Nielsen arc $\alpha$ in
$U_x$ joining $x$ and $y$ such that the image of a homotopy rel $K$
from $\alpha$ to $h(\alpha)$ lies in $V_x$ then  we say that $h$ has
{\em locally constant Nielsen behavior rel $K$ at $x$.} If $h$ has locally constant Nielsen behavior rel $K$ at $x$ for every $x \in \Fix(h)$ then we say that  $h$ has
{\em locally constant Nielsen behavior rel $K$}.
\end{defn}

 If $x \in \Fix(h)$ is either not in $K$ or is an isolated point in $\Fix(h) \cap K$ then $h$ has
locally constant Nielsen behavior rel $K$ at $x$.  One can choose  $U_x \subset V_x$  to be  any contractible neighborhoods   such that $f(U_x) \subset V_x$  and such that $K \cap V_x$ is either empty or $\{x\}$.  

For general $x \in K$, $h$ need not have locally constant Nielsen behavior at $x$  as
the following example illustrates.  Let $h: \R^2 \to \R^2$ be given in
polar co-ordinates by $h(r,\theta) = (r, \theta + f(r))$ where $f(r)$
is a $C^\infty$ function with values between $0$ and $1$ which 
vanishes if and only if $r = 0$ or $r =
1/n,\ = 1,2,3,\dots.$ Choose $\nu_n$ so that $1/(n+1) < \nu_n < 1/n$
and so $f(\nu_n)$ is a non-zero rational multiple of $\pi.$ Then any
point on the circle of radius $\nu_n$ is $h$-periodic with period $>
1$.  Let $K$ be a compact invariant set consisting of the origin
together with a single periodic orbit in the circle of radius $\nu_n$
for each $n$.  There is no Nielsen arc rel $K$ joining the origin to
any other fixed point (the other fixed points must lie on circles of
radius $1/n$).

\begin{defn} \label{defn:nielsen} 
Suppose $h$ has locally constant Nielsen behavior rel $K$ and let $x,y
\in \Fix(h)$. We say that {\em $x$ is Nielsen equivalent to $y$ rel
$K$} if there is a Nielsen arc rel $K$ connecting $x$ to $y$.  The
Nielsen equivalence class of $x \in \Fix(h)$ relative to $K$ is
denoted $\N(h,K,x)$.
\end{defn}

It is immediate that if $h$ has locally constant Nielsen behavior
rel $K$ then each Nielsen class is both open and closed in
$\Fix(h)$.  In particular, if $\Fix(h)$ is compact then  there are only finitely many Nielsen classes
rel $K.$

Having defined relative Nielsen classes in the appropriate context, we now define the    Lefschetz index for a relative Nielsen class.  

\begin{defn} \label{defn:Nielsen-isotopy} 
Suppose $h_t:M \to M$ is an isotopy rel $K$ from $h_0$ to
$h_1$ and suppose $H: M \times I \to M \times I$, given by $H(x,t) =
(h_t(x), t)$, has locally constant Nielsen behavior rel $K \times I.$
Let $\N_0 = \N(h_0,K,x_0)$ be a Nielsen class rel $K$ for $h_0$.  We
define {\em $\N_t$, the Nielsen class rel $K$ determined by $\N_0$ and
$h_t$,} to be $\N_t = \pi( \N \cap (M \times \{t\}))$ where $\N = \N(H,K
\times I,(x_0,0))$ and $\pi : M \times I \to M$ is projection.
\end{defn}

The fact that $H$ has locally constant Nielsen behavior rel $K \times
I$ implies that $h_t$ has locally constant Nielsen behavior rel $K$
for each $t$.  This is because if $\alpha$ is a Nielsen arc in $S
\times I$ joining $(x, t)$ to $(y, t)$ then $\pi \circ \alpha$ is a
Nielsen arc for $h_{t} : S \to S.$ This same observation shows that
$\N_t$ is a Nielsen class rel $K$ for $h_t$ and in fact $\N_t =
\N(h_t, K, x)$ for any $(x,t) \in \N.$

If  $X_0$ is an open and closed subset of $\Fix(h)$ then the Lefschetz 
number $L(X_0,h)$ is well defined.  We recall the definition of $L(X_0,h)$ following \cite{Dold}.   Choose compact $X$ and open $V$ such that $X_0 \subset X \subset V \subset S$ and such that $\Fix(h|V) = X_0$.   We may embed $M$ as a submanifold of $\R^m$
for some $m$ and
consider a retraction $r:W \to M$ where $W$ is a neighborhood of
$S$.  Let $W_V =
r^{-1}(V)$ which is an open neighborhood of $X$ in $\R^m.$ 
There is a canonical ``fundamental'' homology
class $[u] \in H_m(W_V, W_V\setminus X)$.  This class is characterized
by the property that for each $x \in X$ the inclusion $(W_V, W_V
\setminus X) \to (W_V, W_V \setminus \{x\})$ takes $u$ to the standard
generator of $H_m(W_V, W_V \setminus \{x\})$.  The preimage of $0$ under $(id-h \circ r):W_V \to \R^m$ is $X_0$ so  
$$
(id - h \circ r) : (W_V, W_V  \setminus X) \to (\R^m, \R^m \setminus \{0\}).
$$    
The {\em Lefschetz fixed point
index $L(X_0,h)$ of $X_0$} is the integer defined by 
 
\[
(id - h \circ r)_*([u]) = L(X_0,h)[o_m]
\]
where $[o_m] \in H_m(\R^m, \R^m \setminus \{0\})$ is the standard generator.   A proof that  
$L(X_0,h)$ depends only on $X_0$ and $h$ and not on the choice of $X$ and $V$ can be found in \cite{Dold}.

\begin{prop}\label{prop:Nielsen-index}
Suppose $h_t:M \to M$ is an isotopy rel $K$ from $h_0$ to
$h_1$ and suppose $H: M \times I \to M \times I$, given by $H(x,t) =
(h_t(x), t)$, has locally constant Nielsen behavior rel $K \times I.$
Let $\N_0$ be a Nielsen class rel $K$ for $h_0$ and let
$\N_t$ be the Nielsen class determined by $\N_0$ and $h_t$.  Then
$L(\N_t, h_t)$  is independent of $t$.
\end{prop}

\begin{proof}  Let $Y$ be an open neighborhood in $M \times I$
of the compact Nielsen class $\N = \N(H,K \times I,(x_0,0))$ 
where $x_0 \in \N_0$ and assume $Y$ is chosen to be
disjoint from any fixed points of $H$ which are not in $\N$.
We define $V_t = Y \cap (M \times \{t\}).$
It follows from the fact that $\N$ is compact
that for each $t_0$ there is an $\epsilon>0,$ a compact
set $X \subset M$, and an open set $V \subset M$ such that
for each $t \in [t_0 -\epsilon, t_0 +\epsilon]$ we have
$\N_t \subset X \subset V \subset V_t$.

It then follows from the definition  of the Lefschetz index cited above that
$L(N_t,h_t)$ can be defined in terms of $X$ and $V$ for $t \in [t_0 -\epsilon, t_0 +\epsilon]$.  The defining equation becomes $(id - h_t \circ r)_*([u]) = L(\N_t, h_t)   [o_m]$ which depends only on the homotopy class of $h_t$ and is therefore  constant on $[t_0 -\epsilon, t_0 +\epsilon]$.  
Therefore by the compactness of $[0,1]$, $L(N_t,h_t)$    is independent of $t.$
\end{proof}

We now return to the surface case $f : S \to S$.
The following lemma uses smoothness to establish a strong form of locally constant Nielsen behavior relative to $\Fix(f)$.   

\begin{lemma} \label{smoothness}  Suppose that $f \in \Diff_0(S)$, that  $x$ is an accumulation point in $\Fix(f)$ and that $U_{\epsilon}(x)$ is the open $\epsilon$ neighborhood of $x$ in local coordinates around $x$.   Then the following hold for all sufficiently small $\epsilon > 0$ and for  $\alpha_{x,y}$ the   straight line segment  from $x$ to   $y$.   
\begin{enumerate}
\item If $y \in \Fix(f) \cap U_{\epsilon}(x)$ then $f(\alpha_{x,y}) \subset U_{\epsilon}(x)$.
\item There is a neighborhood $W$ of $f$ in $\Diff_0(S)$ so that if $g \in W$, if  $p \notin \Fix(g)$ and if $p,g(p) \in U_{\epsilon}(x)$ then $\alpha_{p,g(p)} \cap \Fix(g) = \emptyset$.
\end{enumerate}
In particular, the straight line homotopy between $\alpha_{x,y}$ and $f(\alpha_{x,y})$ is a homotopy rel $\Fix(f)$ for all  $y \in \Fix(f) \cap U_{\epsilon}(x)$.   
\end{lemma}

\proof   If $y_j \in \Fix(f)$  and $y_j \to x$   then $(y_j - x)/ \|y_j - x\|$ converges to a unit  vector that is fixed by $Df_x$.   If $z_j$ is contained in the interior of $\alpha_{x,y_j}$ then $(z_j - x)/ \|z_j - x\|$ and $(y_j - z_j)/ \|y_j - z_j\|$ converge to this same unit vector.  Thus, for all sufficiently large $j$, the angle between  $\alpha_{x,f(z_j)}$ and $\alpha_{x,z_j}$ 
 is less than one degree and similarly the angle between $\alpha_{y_j,f(z_j)}$ and $\alpha_{y_j,z_j}$   is less than one degree.  In particular $f(\alpha_{x,y_j})$ is contained in the closure of $U_{\epsilon_j}$ where  $\epsilon_j$ is    the distance from $x$ to $y_j$.  This proves (1).   

 If (2) fails  then  there exist  $f_j \to f$,  $p_j \notin \Fix(f_j)$ such that $p_j \to x$  and $q_j \in \Fix(f_j) \cap  \alpha_{p_j,f_j(p_j)}$.  In this case  $(p_j - q_j)/ \|p_j - q_j\|$ converges to a  eigenvector of $Df_x$ with negative eigenvalue.   This contradicts the fact that $f$ preserves orientation and $Df_x$ has an eigenvector  with eigenvalue $1$. 
\endproof 

Motivated by Lemma~\ref{smoothness} we focus on $K$ satisfying the property that each $x \in K$ has a neighborhood that is either disjoint from $\Fix(f)$ or contained in $\Fix(f)$.

\begin{lemma} \label{lem:loc_const_Nielsen}  
Suppose that $f \in \Diff_+(S)$  and that $K$ is a compact invariant set such that $\Fix(f)
\cap K$ is an open (and closed) subset of $K$.  Then $f$ has locally
constant Nielsen behavior rel $K$.  Moreover, if $f_t$ is an isotopy
rel $K$ and $F: S \times I \to S \times I$ is given by $F(x,t) =
(f_t(x), t)$ then each $f_t$ has locally constant Nielsen behavior rel $K$ and 
$F$ has locally constant Nielsen behavior rel $K
\times I.$
\end{lemma}

\begin{proof}   
For the first part, we must show that $f$ has locally
constant Nielsen behavior rel $K$ at $x$ for each $x \in \Fix(f)$ and for this there is no loss in assuming that $x$ is an  accumulation point in $\Fix(f) \cap K$.   We work in local coordinates around $x \in \Fix(f)$, letting $U_{\epsilon}(x)$ be the open $\epsilon$ ball around $x$ and letting $\alpha_{y,z}$ be the straight line segment connecting $y$ to $z$.  
By  Lemma~\ref{smoothness} there exists $\epsilon > 0$ so that, setting $U =  V = U_{\epsilon}(x)$, there is, for each $y \in \Fix(f) \cap U$,  a homotopy rel $\Fix(f)$ from $\alpha_{x,y}$ to  $f(\alpha_{x,y})$   through paths in $V$.    If $\epsilon $ is sufficiently small then $K \cap V \subset \Fix(f) \cap V$ and so the homotopy is also rel $K$.  

We proceed to the argument for $F$.   Suppose that $0 \le t \le 1$ and that $(x,t) \in \Fix(F)$.    Note that $K \cap \Fix(f_s)$ is independent of $s$  since the isotopy is rel $K$.  Thus  each $f_s$ has locally constant Nielsen behavior rel $K$.    Moreover, there is no loss in assuming that   $x$ is an accumulation point in $\Fix(f_t) \cap K$.    
For $\gamma > 0$ let  $J_{\gamma}(t)$ be the $\gamma$ neighborhood of $t$ in  $I$.  
 Choose  $\epsilon_t$ as in the last sentence of the preceding paragraph  with $f$ replaced by $f_t$  and choose $\delta_t$ so that $f_t(U_{\delta_t}(x)) \subset U_{\epsilon_t/2}(x)$.   Denote $U_{\delta_t}(x) \times J_{\gamma}(t)$ by $U$ and   $U_{\epsilon_t}(x) \times J_{\gamma}(t)$ by $V$.    Thus $U \subset V$ and if $\gamma$ is sufficiently small then  $F(U) \subset V$.   If  $(y,s) \in \Fix(F) \cap U$ let $\beta$ be the path in $U$ connecting $x$ to $y$ that is 
 the concatenation
of the ``vertical'' arc $\{x\} \times [t,s]$ with the  ``horizontal'' arc $\alpha_{x, y} \subset U \times \{s\}$.  To prove that $F(\beta)$ is homotopic to $\beta$ rel $K \times I$ through paths in $V$ it suffices to show that the straight line homotopy from $f_s(\alpha_{x, y})$  to  $\alpha_{x, y}$ in $U_{\epsilon_t}(x)$ is a homotopy relative to $K$.  This follows from Lemma~\ref{smoothness}(2) if $\gamma$ is sufficiently small.
\end{proof}

We conclude this section with a variant of Proposition~\ref{prop:lift}.

\begin{lemma} \label{applying the pivot lemma} Suppose that $f \in 
\Homeo_0(S)$, that $K$ is a compact $f$-invariant set and $D \subset S$ is 
a disk whose boundary is disjoint from $K$ and whose interior intersects 
$K$ in a non-empty set that is disjoint from $\Fix(f)$.  Suppose further 
that there is a lift $\ti f : \ti S \to \ti S$, a component $\ti D_0$ of the full pre-image of $D$ in $\ti S$ and an isotopy rel  $\ti K$, the complete pre-image of $K,$ from $\ti f$ to a homeomorphism that preserves $\ti D_0$.   
 Then there
exists $\ti x \in \Fix(\ti f)$ whose projected image  $x \in \Fix(f)$   is not Nielsen equivalent rel $K$ to any element
of $K$.  If  $\ti f$ is the identity lift then $x \in \FC(f)$.
\end{lemma}

\proof The compact set $C= \ti K \cap \ti D$ is $\ti f$-invariant. For any 
$\ti y \in \ti K \cap \Fix(\ti f)$, our hypotheses imply that $\ti y$ is 
disjoint from $\ti D$ and that if $R_{\ti y} \subset \ti S$  is a properly embedded ray $R_{\ti y} \subset \ti S$ initiating at 
$\ti y$  that is disjoint from $\ti D$ then  $\ti f(R_{\ti y})$ is properly homotopic to $R_{\ti y}$ rel $C$ and rel endpoint.  
 By Proposition~5.3 of 
\cite{FLP}, there exists $\ti x \in \Fix(\ti f)$ for which 
there 
does not exist a properly embedded ray $R_{\ti x} \subset \ti S$ initiating at $\ti x$ such that $\ti f(R_{\ti x})$ is properly homotopic to $R_{\ti x}$ rel $C$ and rel endpoint.   It follows that 
there does not exist 
 an arc $\ti \alpha$ connecting $\ti x$ to   $\ti y$
 such 
that $\ti f(\alpha)$ is homotopic to $\alpha$ rel $\ti C$ and rel 
endpoints.    The projected image $x \in \Fix(f)$ of $\ti x$ is therefore not Nielsen 
equivalent rel $K$ to any $y \in K$.
\endproof

\section{Reducing Curves} \label{rc}

We assume throughout this section that $f_1,\ldots,f_m$ generate  an abelian subgroup $\F \subset \Diff_0(S)$.   

 For each $x \in S$, define  $A_x \subset
\{1,\ldots,m\}$ by $i \in
A_x$ if and only if  $x \in \Fix(f_i)$.  The cardinality of $A_x$ is called the {\em multiplicity} of $x$.  A   set $K$ has {\em locally constant multiplicity} if it satisfies the equivalent conditions of the following lemma.

\begin{lemma}  The following are equivalent for any   set $K$.
\begin{enumerate}
\item Each $x \in K$ has
a neighborhood $U$ such that the multiplicity of $y \in K \cap U$ is
independent of $y$. 
\item Each $x \in K$ has a neighborhood $U$ such
that $A_y$ is independent of $y \in K \cap U$.  
\item  $\Fix(f_i) \cap K$ is an open subset of $K$ for all $i$,
\end{enumerate}
\end{lemma}

\begin{proof}
Let $\cP$ denote the power set of $\{1,\ldots,m\}$.  
Property (2) asserts that the function $\cA : K \to \cP$ given
by $\cA(x)  = A_x$ is locally constant. Likewise, property
(1) asserts that the function $card_{\cA} : K \to \Z$ is locally
constant where $card_{\cA}(x)$ is the cardinality of $A_x$.  
Finally (3) asserts that for each $i$ the function $\psi_i : K
\to \{0,1\}$ is locally constant, where $\psi_i(x) = 1$ if
$x \in \Fix(f_i)$ and $0$ otherwise.

To see that (2) implies (3) observe that $\psi_i = \phi_i \circ \cA$
where $\phi_i : \cP \to \{0,1\}$ is the function defined by
$\phi_i(A) = 1$ if $i \in A$ and $0$ otherwise. 

To see that (3) implies (1) observe that the function 
\[
card_{\cA}(x) = \sum_{i=1}^m \psi_i(x).
\]

To see that (1) implies (2) note the fact each $\Fix(f_i)$ is a closed
set implies the function $\cA$ is semicontinuous in the sense that if
$\lim x_n = x$ and $\cA(x_n) = A_0$ is independent of $n$ then $A_0
\subset A_x.$ If in addition $card_{\cA}$ is locally constant it
follows that the cardinality of $A_0$ equals the cardinality of $A_x$
so $A_0 = A_x.$ This implies that the function $\cA$ is continuous or
equivalently locally constant.
\end{proof}

\begin{lemma}\label{lem:K_exists}
Suppose that $Y$ and $C$ are compact
disjoint $\F$-invariant subsets of $S$ and that $Y \subset \cup_{i=1}^m \Fix(f_i)$.  Then there is a compact
$\F$-invariant subset $K \subset Y$ with locally constant multiplicity
such that for all $i$, each $y \in Y \cap \Fix(f_i)$ is $f_i$-Nielsen equivalent
rel $K \cup C$ to some $x \in K$.
\end{lemma}

\proof The proof is by induction on the maximal multiplicity $M(Y)$ of an element of $Y$.  If $M(Y)= 1$ then $K = Y$ satisfies the conclusions of the lemma.  We may therefore assume that  the lemma holds  for   $Y'$ with $M(Y') < M(Y)$.  Denote $M(Y)$ by $M$ and let $K_M$ be the set of  points in $Y$ with multiplicity $M$. 

If $x \in K_M$ is isolated in $K_M$, let $U_x$ be an open neighborhood of $x$ that intersects $K_M$ in $\{x\}$.  Otherwise, $x$ is an accumulation point in $K_M$ and hence an accumulation point in $\Fix(f_i) \cap K$ for each $i \in A_x$.  In this case we let $U_x$ be an open  neighborhood of $x$ such that $U_x \cap \Fix(f_i) = \emptyset$ if $i \notin A_x$ and such that for each $i \in A_x$ and each $y \in U_x \cap \Fix(f_i)$ there is a   path $\beta_{x,y}$ connecting $x$ to $y$ in $U_x$ and a homotopy rel $\Fix(f_i)$ from $\beta_{x,y}$ to $f_i(\beta_{x,y})$ through paths in $U_x$.  Such a neighborhood exists by  Lemma~\ref{smoothness}.  

 The $U_x$'s cover $K_M$ and so have a finite subcover $\{U_{x_1},\ldots,U_{x_k}\}$.  Define $W$ to be  the $\F$-orbit of $\cup_{j=1}^k U_{x_j}$,    $Y'$ to  be the intersection of $Y$ with the complement of $W$ and $C'$ to be $C \cup K_M$.   The inductive hypothesis applies to $Y'$ and $C'$ to produce $K' \subset Y'$;   let $K = K' \cup K_M$.   

$K$ has locally constant multiplicity because both $K'$ and $K_M$ do and are disjoint.  It suffices to show that each $y \in Y \cap \Fix(f_i)$  is $f_i$-Nielsen equivalent rel $K \cup C$  to some $x \in K$. This follows from the inductive hypothesis if $y \in Y'$.  We may therefore assume that $y \in W$.  Suppose at first that $y \in U_{x_j}$ for some $j$.  Since $\Fix(f_i) \cap U_{x_j} \supset  K_M \cap U_{x_j} = (K \cup C) \cap U_{x_j}$   the homotopy from $\beta_{x_j,y}$ to $f_i(\beta_{x_j,y})$  is rel $K \cup C$.  Thus $y$ is $f_i$-Nielsen equivalent to $x_j$ rel $K \cup C$.  For the general case we must consider $f(y)$ where $y \in U_{x_j}$ and $f \in \F$.   Since $K,C$ and $\Fix(f_i)$ are all $f$-invariant and since $f$ commutes with $f_i$, the path $\gamma = f(\beta_{x_j,y})$ connecting $f(x_j)$ to $f(y)$ is homotopic rel $K \cup C$ to $f_i(\gamma)$ showing that  $f(y)$ is $f_i$-Nielsen equivalent to $f(x_j)$ rel $K \cup C$.  
\endproof

\begin{defn}  \label{def:reducing}A {\em reducing set}  for $\{f_1,\ldots,f_m\}$ relative to a compact $\F$-invariant set $K$ is a finite collection $\cR$ of disjoint simple closed curves in $S$ with the following properties.  
\begin{enumerate} 
\item $\cR \cap K = \emptyset$ and $\cR$ is $\F$-invariant up to isotopy rel $K$.
\item Each contractible component of  $S \setminus \cR$ contains at least two elements of $K$.
\item  For $i =1,\ldots,m$ the following holds:
\begin{description}
\item[(3-i)] If a component   $S_l$ of $S \setminus \cR$   intersects $K$ in an infinite set and if $S_l \cap K$ contains a fixed point of $f_i$, then $f_i$ is isotopic rel $K$ to a homeomorphism that is the identity on $S_l$.
\end{description}
\end{enumerate}
\end{defn}

Note that if there is a reducing set for $\{f_1,\ldots,f_m\}$ then each $K \cap \Fix(f_i)$ is an open subset of $K$ and so $K$ has locally constant multiplicity.

\begin{defn}  If $\N_i$ is a union of Nielsen classes for $f_i$ then a compact set $K$ is {\em $\N_i$-complete}  if each element of $\N_i$ is Nielsen equivalent rel $K$ to an element of $K$. 
\end{defn}

\begin{proposition}\label{prop:reducing}
Suppose that $\N_i$ is a union of Nielsen classes for $f_i$ and that
$K \subset \cup_{i=1}^m\ \N_i$ is compact, is $\F$-invariant, has locally
constant multiplicity and is $\N_i$-complete for each $1 \le i \le m$.
Then there is a reducing set $\cR$ for $\{f_1,\ldots,f_m\}$ rel $K$.
Moreover if $\cR'$ is a finite set of disjoint simple closed curves satisfying (1) and (2) of Definition~\ref{def:reducing} then we may choose $\cR$ to contain $\cR'$.
\end{proposition}

\proof  We show below that there is a finite collection  $\cR_m$ of disjoint simple closed curves containing $\cR'$ and satisfying (1), (2) and (3-m) of Definition~\ref{def:reducing}.  This proves the proposition in the case that $m=1$.
Assume by induction on $m$ that the proposition holds for $\{f_1,\ldots,f_{m-1}\}$. To produce the desired reducing set $\cR$   for $\{f_1,\ldots,f_m\}$ rel $K$,  apply the inductive statement to $\{f_1,\ldots,f_{m-1}\}$ with  $R'$   replaced by $\cR_m$.

We now prove the existence of $\cR_m$, denoting $f_m$, $\cR_m$ and $\N_m$ simply by $f$, $\cR$ and $\N$.  The proof is based on arguments from section 2 of \cite{FHP}.  Both $B= \Fix(f) \cap K$ and  $C = K \setminus B$  are $\F$-invariant and are compact by the assumption that $K$ has locally constant multiplicity.  If $B$ is finite then simply take $\cR = \cR'$. We assume now that $B$ is infinite.

Define   $\WF$ to be the set of compact  subsurfaces   $W \subset S$ such that:
\begin{description}
\item [($\W$-1)] $\partial W$  has finitely many components, each of which is disjoint from $K$ and is either disjoint from $\cR'$ or is isotopic rel $K$ to an element of $\cR'$.  Moreover, each contractible   component of $S \setminus \partial W$ contains at least two points in $K$.
\item [($\W$-2)]  $W$ contains all but finitely many points of $B$ and  every component of $W$ intersects $B$ in an infinite set.
\item [($\W$-3)] $f$ is isotopic rel $K$ to a homeomorphism that is the identity on $W$.
\end{description}

 Lemma 4.1 of \cite{han:commuting} and the isotopy extension theorem (see the proof of Theorem~1.2 of \cite{fh:periodic} where $\WF$ is denoted $\W(\langle f \rangle,B, K \setminus B)$; see also Lemma~\ref{smoothness} of this paper) imply  that $\WF \ne \emptyset$.  (This uses the fact that $B$ is compact and is open in $K$.)
There is a partial order on $\WF$ defined by  $W_1 < W_2$ if and only if $W_1$ is isotopic rel $K$ to a subsurface of $W_2$ but is not isotopic rel $K$   to $W_2$.   

The next step in the proof is to show   that $\WF$ has  maximal elements.  This is essentially Lemma~2.7 of \cite{FHP}; the main difference is that instead of assuming that $B$ contains each element of $\N$ we are assuming that every element of $\N$ is Nielsen equivalent rel $K$ to an element of $B$.  To compensate for this, we quote Lemma~\ref{applying the pivot lemma} instead of the Brouwer lemma to produce a fixed point in a certain universal cover.  We make this precise as follows.

 If $W_1 \in \WF$ is not contained in a maximal element of $\WF$ then there is an infinite increasing sequence  $W_1 \subset W_2 \subset \dots$ of elements of $\WF$ that are non-isotopic rel $K$.  We may assume that $W_l \cap B$ and the number of components of $W_l$ are independent of $l$.   The number of components of $S \setminus W_l$ is unbounded.   We may therefore choose $l$ so that some complementary component of $W_l$ is a disk $D$ that is disjoint from $B$. 

There is a homeomorphism  $\phi : S \to S$ that is isotopic to $f$ rel $K$ and such that
$\phi|_{W_l}$ is the identity.  Choose a component $\ti W_l^*$ of the full pre-image $\ti W_l$ of $W_l$ in $\ti S$ and let $\ti \phi$ be the lift of $\phi$ that restricts to the identity on $\ti W_l^*$.   The isotopy rel $K$ from $\phi$ to $f$ lifts to an isotopy rel $\ti K$ from $\ti \phi$ to some lift $\ti f$ of $f$.   Choose $y \in B \cap W_l$ and a lift $\ti y \in \ti W_l^*$.  Then $\ti y \in \Fix(\ti \phi)$ and, since $\ti \phi$ is isotopic to $\ti f$ rel $\ti K$, $\ti y \in \Fix(\ti f)$.  Thus every element of $\Fix(\ti f)$ projects into $\N \subset K$. 

 On the other hand,  some component $\ti D^*$ of the the full pre-image of $D$ in $\ti S$   is a complementary component of $\ti W_l^*$.   Lemma~\ref{applying the pivot lemma} implies that there exists $\ti x \in \Fix(\ti f)$ whose projected image $x \in \Fix(f)$ is not Nielsen equivalent rel $K$ to an element of $K$.   This contradiction completes the proof that  $\WF$ has  maximal elements.

The hypotheses of Lemma~2.7 and Lemma~2.8 of \cite{FHP} are slightly different than those of this proposition but their proofs carry over without change to this context.  Together they imply that up to isotopy rel $K$ there is a unique maximal element $W$ of $\WF$ and that one may choose $W$ so that $\partial W$ is disjoint from $\cR'$.   Let $\cR_m = \partial W \cup \cR'$. Properties ~\ref{def:reducing}(2) and ~\ref{def:reducing}$(3-m)$ follow from the fact that $W \in \WF$.  For ~\ref{def:reducing}(1), note that $g(W)$ is also a maximal element of $\WF$  for each $g \in \F$ because  $gfg^{-1} = f$.  The  uniqueness of $W$ up to isotopy rel $K$ implies that $W$, and hence $\cR_m$, is $g$-invariant up to isotopy rel $K$.
\endproof

\begin{cor} \label{reducing sets exist} 
Suppose that $\{f_1,\ldots,f_m\}$ generate an abelian subgroup of
diffeomorphisms of $S$.    Then there is a compact subset
$K$ of $\cup_{i=1}^m \FC(f_i)$ that has locally constant Nielsen
behavior for each $f_i$ and that is $\FC(f_i)$-complete for each $i$ and
there is a reducing set $\cR$ for $\{f_1,\ldots,f_m\}$ relative to $K$.
\end{cor}
 
\begin{proof}
Ley $Y = \cup_{i=1}^m \FC(f_i)$.   The compact 
$\F$-invariant subset $K \subset Y$  produced by Lemma~\ref{lem:K_exists} applied to $Y$ 
and $C = \emptyset$ has locally constant multiplicity and is $\FC(f_i)$-complete
for each $i.$  The former implies by Lemma~\ref{lem:loc_const_Nielsen}   that each $f_i$ has
locally constant Nielsen behavior rel $K$.
The existence of a reducing set $\cR$ for  $\{f_1,\ldots,f_m\}$ relative to $K$
  now follows from Proposition~\ref{prop:reducing}.
\end{proof}

\section{Thurston Normal Form}

Our definition of a reducing set $\cR$ for $\{f_1,\ldots,f_m\}$ rel $K$ places no restrictions on the components $Y$ of $S \setminus \cR$ that intersect $K$ in a finite set.  The natural restriction is that the isotopy class of $Y$ relative to $K$ determined by $f_l$ should be  irreducible in the sense of Thurston normal form (TNF) for each $f_l$.  In this section we recall the appropriate definitions and a result about Thurston normal form for abelian subgroups.   The actual restrictions will be made in the next section.  Further details about Thurston normal form can be found, for example, in \cite{Th}, \cite{FLP}, \cite{HT}. 
 
Let $S$ be a finitely punctured surface of negative Euler characteristic, perhaps with boundary and let $\MCG(S)$ be the  mapping class group of $S$.  

Suppose that $\cR_T$ is a  finite set of disjoint simple closed curves on $S$, all of whose complementary components $S_i$ have negative Euler characteristic and that $f : S\to S$ is a homeomorphism that setwise preserves $\cR_T$. For each  $S_i$  there exists a minimal $n_i \ge 1$ such that $f^{n_i}$ preserves $S_i$.  Let $\bar S_i$ be the completion of $S_i$ obtained by compactifying each end of $S_i$ that corresponds to an element of $\cR_T$ with a boundary circle.  Then $f^{n_i}|S_i$ determines an element $\alpha_i \in \MCG(\bar S_i)$ that depends only on the isotopy classes of the elements of $\cR_T$ and the isotopy class $\alpha \in \MCG(S)$ determined by $f$.  We say that $\cR_T$ is {\em a TNF-reducing set for $f$ (or $\alpha$)} if each $\alpha_i$ is either represented by a periodic homeomorphism or by a pseudo-Anosov homeomorphism. 

Let $\cS$ be the set of isotopy classes of simple, closed, essential,
non-peripheral curves on $S$.  Then $\MCG(S)$ acts on $\cS$ and for
each element $\alpha$ of $\MCG(S)$ we define $\Per(\alpha)$ to be the
set of elements of $\cS$ on which $\alpha$ acts periodically.  The
intersection number $\gamma\cdot \delta$ of two elements
$\gamma,\delta \in \cS$ is the minimum, over all representatives $c$
of $\gamma$ and $d$ of $\delta$, of the cardinality of $c \cap d$.
Following \cite{HT}, we define $\Per_0(\alpha)$ to be those elements
$\gamma$ of $\Per(\alpha)$ such that $\gamma\cdot \delta = 0$ for all
$\delta \in \Per(\alpha)$.

\begin{lemma} $P_0(\alpha)$ is a TNF-reducing set for $\alpha$.
\end{lemma}

\proof   We assume notation as above.  By construction, $\Per_0(\alpha_i) = \emptyset$ for all $\alpha_i$.  Lemma 2.2 of \cite{HT} (see also Lemma~2.12 of \cite{FHP}) asserts that if $\alpha_i$ is not represented by a periodic homeomorphism then $\Per(\alpha_i) = \emptyset$.  It
then follows from the Thurston normal form theorem that $\alpha_i$ is represented by a pseudo-Anosov homeomorphism. 
\endproof 

\begin{cor}\label{cor: ht}  If $A$ is an abelian subgroup of $\MCG(S)$ then $\bigcup_{\alpha \in A} \Per_0(\alpha)$ is a finite $A$-invariant subset of $\cS$ that is represented by disjoint simple closed curves and is a  TNF-reducing set for each $\alpha \in A$.
\end{cor}

\proof If $\alpha, \beta \in A$ then $\beta(\Per_0(\alpha)) = \Per_0(\alpha)$
since  the action of $\beta$ on $\cS$ preserves $\Per(\alpha)$ and intersection numbers.  It follows that  $\Per_0(\alpha) \subset \Per(\beta)$ and hence that $\gamma \cdot \delta=0$ for all $\gamma \in \Per_0(\alpha)$ and $\delta \in \Per_0(\beta)$.   This proves that  $\bigcup_{\alpha \in A} \Per_0(\alpha)$ is a finite subset of $\cS$ and since it contains a TNF-reducing set for each $\alpha \in A$  it is a TNF-reducing set for each $\alpha \in A$.
\endproof 

\begin{remark}  
It is easy to see that $\Per_0(\alpha)$ is the minimal TNF-reducing
set for $\alpha$ in the sense that every element of any TNF-reducing
set for $\alpha$ is an element of $\Per_0(\alpha)$.  It follows that
$\bigcup_{\alpha \in A} \Per_0(\alpha)$ is a 
minimal subset of $\cS$ that is a TNF-reducing set for every $\alpha \in A$.
This minimal subset has a canonical characterization which we now
describe.  Define $\Per(A)$ to be those elements of $\cS$ which
are in $\Per(\alpha)$ for all $\alpha \in A$ 
and $\Per_0(A)$ to be those elements $\gamma \in \Per(A)$ such that
$\gamma\cdot \delta = 0$ for all $\delta \in \Per(A)$. 
Since for all $\alpha,\beta \in A$ the set $\Per_0(\alpha)$ is
finite and $\beta(\Per_0(\alpha)) = \Per_0(\alpha)$ it is
clear that $\Per_0(\alpha) \subset \Per_0(A)$ for
every $\alpha \in A.$  Therefore $\Per_0(A)$ 
is a TNF-reducing set for every $\alpha \in A$.
It is straightforward to show $\Per_0(A)$
is also minimal so $\bigcup_{\alpha \in A} \Per_0(\alpha) =
\Per_0(A)$.

\end{remark}

\section{Proof of Theorem~\ref{thm:main}}

We assume throughout this section that $\F$ is an abelian subgroup of
$\Diff^+(S)_0$ with rank $k$.

\begin{defn}
A subset $\Gamma$ of $\F$  will be called
{\em $k$-spanning} if it has cardinality at least $k$ and if every 
subset $\Gamma_0$ of cardinality $k$ generates a finite index subgroup
of $\F.$  Equivalently $\Gamma$ is $k$-spanning if every 
subset $\Gamma_0$ of cardinality $k$ is a basis of the vector
space $\R^k$ determined by the identification of $\F$ with $\Z^k$.
\end{defn}

\begin{lemma} \label{lem:basis} 
Given $m \ge k$ there is a subset $\Gamma \subset \Z^k$ of cardinality
$m$ such that $\Gamma$ is $k$-spanning.
\end{lemma}
\begin{proof}
The desired property is equivalent to the requirement that 
$\Gamma$ contains a set of generators and that any
subset of $\Gamma_0 \subset \Gamma$ of cardinality $k$ is 
linearly independent in $\R^k.$
We induct on $m$.  For $m = k$ take $\Gamma$ a set of generators of
$\Z^k.$  Given a $\Gamma_m$ of size $m$ with the desired property
let $X \subset \R^k$ be the
union of the spans of all subsets of $\Gamma_m$ with cardinality $k-1$.
Then $X$ is a union of a finite number of $(k-1)$-dimensional subspaces.
Let $v$ be an element of $\Z^k$ which is not in $X$ and define
$\Gamma_{m+1} := \Gamma_m \cup \{v\}.$  If $\Gamma_0 = \{v, v_1, \dots
v_{k-1}\}$ then $v$ is not a linear combination of 
$\{v_1, \dots v_{k-1}\}$ so $\Gamma_0$ is linearly independent and
$\Gamma_{m+1}$ has the desired property.
\end{proof}

We establish notation for the remainder of this section as follows.  $\F_m$ is a finite index subgroup of $\F$ generated by a $k$-spanning set $\{f_1, \dots, f_m\} \subset \F$.  $K$ is a compact set that has locally constant Nielsen behavior with respect to   $f_i$ for $1 \le i \le m$.   $\cR$  is a reducing set for $\{f_1, \dots, f_m\}$ relative to  $K$.    

Each  $f \in \F_m$ is isotopic rel $K$ to a homeomorphism $F$ that preserves $\cR$ and so preserves each complementary component  of $\cR$.  Suppose that $Y$ is such a component, $K_Y = K \cap Y$ is finite and that $\bar Y$ is the completion of $Y$ obtained by compactifying each end of $Y$ that corresponds to an element of $\cR$ with a boundary circle.      The restriction $F|Y$ determines an element $[f|\bar Y]$ of the mapping class group $\G_{\bar Y}$ rel $K_Y$  that depends only on $f$ and not the choice of $F$.   The assignment $f \mapsto [f|\bar Y]$ defines a homomorphism  $\phi_Y : \F_m \to \G_{\bar Y}$.  

\begin{defn}  Assume notation as above.  $\cR$ is {\em TNF-complete} for  $\{f_1, \dots, f_m\}$  if $[f_i|\bar Y]$ is either periodic or pseudo-Anosov  for each $f_i$ and for each  component $Y$  with finite $K_Y$. If $\cR$ is not TNF-complete then it is {\em TNF-incomplete}.
\end{defn}

\begin{lemma}  \label{reducing set}  Any reducing set $\cR$ for $\{f_1, \dots, f_m\}$ relative to  $K$   extends to a TNF-complete reducing set for $\{f_1, \dots, f_m\}$ relative to  $K$.
\end{lemma}

\proof Suppose that $Y$ is a complementary component of $\cR$ and that
$K_Y$ is finite.  Define $\phi_Y : \F_m \to \G_{\bar Y}$ as above.  Since
$\phi_Y(\F_m)$ is an abelian subgroup of $\G_{\bar Y}$, by
Corollary~\ref{cor: ht} there is a finite set $R_Y$ of disjoint
non-parallel, non-peripheral, simple closed curves in $Y \setminus
K_Y$ so that each element $\phi_Y(\F_m)$ is represented by a
homeomorphism of $Y$ whose restriction to each component of the
complement of $R_Y$ is either periodic or pseudo-Anosov.   Taking the union of $\cR$ with $R_Y$ for
all such choices of $Y$ produces the desired TNF-complete reducing
set.  \endproof

  We assume throughout that $\cR$ is TNF-complete  and that no proper subset of $\cR$ is   a TNF-complete reducing set.    

Let $\E$ be the subset of $\cR$ consisting of those elements that are
essential in $S.$
Since $\chi(S) <0 $ equals the sum of the Euler characteristics of
the components of the complement of $\E$ in $S$ there must be
such a component with negative Euler characteristic.  In fact
if $\{N_i\}_{i=1}^r$ is the set of closures of those components of 
$S \setminus \E$ which have negative Euler characteristic then
$1 \le r \le -\chi(S).$  The remaining components of $S\setminus \E$ are annuli and we let $\{A_j\}_{i=1}^s$ be those annuli that share a boundary component with some $N_i$.  Three distinct $A_j$ cannot be parallel in $S$ so $s$ is bounded by $2-\chi(S)$, which is twice the genus of $S$ and hence twice the maximum cardinality of a set of disjoint essential non-parallel curves in $S$.  Let $S_- = (\bigcup_{i=1}^r N_i) \cup (\bigcup_{j=1}^s A_j)$,\  $S_0 = S \setminus S_-$,\   $K_- = K \cap S_-$ and $K_0 = K \cap S_0$.  The number $|S_-|$ of subsurfaces $N_i$ and $A_j$ that compose $S_-$ is bounded above by $2- 2 \chi(S)$.

Before beginning the formal proof of  Theorem~\ref{thm:main}, we prove it  in a pair of special cases.  

\begin{lemma}\label{must intersect} If $K$ is $\FC(f_l)$-complete  then $\FC(f_l) \cap K_- \ne \emptyset$. 
\end{lemma}

\proof   To simplify notation we write $f_l$ simply as $f$.  We assume that $\FC(f) \cap K_- = \emptyset$ and argue to a contradiction.  Let  $K_0 = K \setminus K_-$.   The first step in the proof is to choose  $g$  isotopic to $f$ and satisfying the following properties.  
\begin{enumerate}
\item $f$ agrees  with $g$ on the complement of $\bigcup A_j$.
\item The restriction of $g$ to each $A_j$ has a periodic
orbit $P_j$.
\item $\Fix(f) = \Fix(g)$ and $\FC(f) = \FC(g)$.
\item \label{still complete} If $K' = K_0 \bigcup P_j$
then every point of $\FC(g)$ is Nielsen equivalent rel $K'$ to
a point of $K'$.
\end{enumerate}

Choose
an $f$ recurrent point $p_j$ in $K \cap \Int(A_j)$.  If $p_j$ is periodic,
 let $P_j = orb(p_j)$   and let $f|A_j =g|A_j$. Note that if $p_j$ is fixed then it is not in $\FC(g)$.   
Otherwise, perturb
$f$ in a neighborhood $W_j$ of $p_j$ in $A_j$ that is disjoint from its $f$-image so that $p_j$ becomes periodic with period
greater than $1$ and let $P_j = orb(p_j).$   Properties (1) through (3) are satisfied.

The fact that we can achieve Property~\ref{still complete} requires justification.   For each $x \in \FC(f)$ choose  neighborhoods $U_x \subset V_x$ as in the definition of locally constant Nielsen behavior rel $K$. Since $\FC(f)$ is compact there exist $x_1,\ldots,x_m$ such that $ \FC(f) \subset \cup_{l=1}^m U_{x_l}$.    Let $\beta_l$ be an $f$-Nielsen arc rel $K$ connecting $x_l$ to some element of $\FC(f) \cap K$.  The existence of $\beta_l$ follows from our assumption that   $K$ is $\FC(f_l)$-complete. If the $W_j$'s are sufficiently small then $\bigcup W_j$ is  disjoint from each $V_{x_l}$ and the homotopy from $f(\beta_l)$ to $\beta_l$  is rel $K \cup (\bigcup W_j)$.  For $x \in U_{x_l} \cap \FC(f)$, let $\beta_x$ be an $f$-Nielsen arc that is the concatenation of an $f$-Nielsen arc in $U_{x_l}$ connecting $x$ to $x_l$ with $\beta_l$.   Then $\beta_x$ is  an $f$-Nielsen arc rel $K \cup  (\bigcup W_j)$ and so is a $g$-Nielsen arc rel $K'$.

Let $\cR_{T}$ be a minimal set of reducing curves (with respect to Thurston normal form) for the element $[g|S_-]$ determined by $g$ in the mapping class group of $S_-$ relative to the finite set $\cup_j P_j$.  After an isotopy we may assume that each essential element of $\cR_T$ is a component of some $\partial N_i$. There is an isotopy rel $K'$ from $g$ to a homeomorphism $G$ that preserves $S_-$,  preserves $\cR_T$ and that is the identity on each $N_i$.   

 For each $A_j$, let $X_j$ be the complement  in $A_j$ of the disks bounded by the inessential elements, if any, of $\cR_T$ that are contained in $A_j$.

We claim that $\cup_i\partial N_i \subset \cR_T$.  
If this fails, then there exists adjacent $N_i$ and $X_j$  such that  $G|(N_i \cup X_j)$ is the identity.   
If $A_j = X_j $  then  $P_j =\{p_j\}$ is a single point $p_j \in \FC(G) \cap K_-$.  Since $G$ is isotopic to $g$ rel $P_j$, $p_j \in \FC(g) = \FC(f)$ in contradiction to the assumption that $\FC(f) \cap K_- = \emptyset$.  
If $A_j \ne X_j $  then  $A_j \setminus X_j$ is a single disk  $D_j$   that contains $P_j$ contradicting property (4) above and    Lemma~\ref{applying the pivot lemma} (applied with $K'$ replacing $K$) which implies the existence of  an element of $\FC(g)$ that is not $g$-Nielsen equivalent rel $K'$ to any element of $K'$.  This completes the proof that $\cup_i\partial N_i \subset \cR_T$.  

    It is well known that    $N_i$ is a complete Nielsen class for $G|S_-$ relative to  $\cup P_j$ and that the action of $G$ on $\pi_1(X_j, \partial X_j)$ does not fix any element determined by an arc in $X_j$ connecting distinct  components of $\partial X_j$.   It follows that $N_i$ is a complete Nielsen class $\mu'$ for $G$ relative to $K'$.    Indeed, if there were a Nielsen path connecting $x \in N_i$ to $y \in S_0$ then its intersection with some $X_j$ would contain an arc with endpoints on distinct components of $A_j$ whose  relative homotopy class would be fixed by the action of $G$.

 On the other hand,  $\mu' \subset \FC(G)$  has Lefschetz index equal to $\chi(N_i) \ne 0$.  Since $\Fix(f) = \Fix(g)$ and since $K'$ is the union of $K_0$ with a finite set, $K' \cap \Fix(g)$ is an open subset of $K'$.  Lemma~\ref{lem:loc_const_Nielsen} and Proposition~\ref{prop:Nielsen-index} imply that $\mu$, the $g$-Nielsen class rel $K'$ determined by $\mu'$ and the isotopy rel $K'$ from $g$ to $G$, is   non-empty.  Since $\mu' \subset \FC(G)$, $\mu \subset \FC(g)$.  By  Property~\ref{still complete},  $\mu$ contains an element   $y$ that is contained in $K'$ and so in not contained in $N_i$.   Since $y \in \Fix(g) \cap K'$  and $G$ is isotopic to $g$ rel $K'$, it follows that $y$ is an element of  $\mu'$.  This contradiction completes the proof.  
\endproof

For each $M = N_i$ or $M=A_j$ there is a component $X_M$ of $S \setminus \cR$ that is either equal to $M$ or is obtained from $M$ by removing disks  whose boundaries  are  elements of $\cR \setminus \E$.

\begin{lemma} \label{many fixed points}
Suppose that $M = N_i$ or  $M=A_j$ and that $\FC(f_l) \cap K \cap X_M  \ne \emptyset$ for  at least $k$ distinct values of $l$.   Then $\Fix_c(\F) \ne \emptyset$.
\end{lemma}

\proof       For simplicity we denote $X_M$ by $X$.  
Let $\F'$ be  the finite index subgroup of $\F$ generated by the $f_l$'s satisfying $\FC(f_l) \cap K \cap X  \ne \emptyset$.    If  each $f_l|X$ is isotopic rel $K$ to a homeomorphism whose restriction to $X$ is the identity then $K \cap X  \subset \FC(f_l)$ for each $l$ and $K \cap X$ is a non-empty subset of  $\FC(\F'))$.  Lemma~\ref{lem:contractible} then completes the proof.
  
We may therefore assume that some $f_l$, which we denote simply as $f$,  is not isotopic rel $K$ to a homeomorphism whose restriction to $X$ is the identity.   It follows that $X \cap K$ is finite.  Since $\cR$ is TNF-complete, $f|X$ determines an irreducible relative mapping class.   Thus $f$ is isotopic rel $K$ to a homeomorphism $F$ whose restriction to $X$ is either periodic with period greater than one or is pseudo-Anosov.    At least two components of $\partial X$ are contained in $\E$ and so are $F$-invariant.  There is also at least one $F$-invariant puncture in $X$.   It follows that $F|X$ is not  periodic with period greater than one.   

It remains to consider the case that $F|X$  is pseudo-Anosov.  Each element of $\F'$ is isotopic to a homeomorphism that preserves $X$; restricting   this homeomorphism to $X$ defines a homomorphism $\phi : \F' \to \G$  where  $\G$ is the subgroup  of the mapping class group of $X$ rel $X \cap K$ consisting of elements that commute with the pseudo-Anosov mapping class $\phi(f)$.  It is well known (see for example Lemma~2.10 of \cite{FHP})  that $\G$ is virtually cyclic.  The subgroup $\F_0 =  \phi^{-1}(\langle \phi(f) \rangle)$
has finite index in $\F$.   Moreover,   $K \cap \FC(f) \cap X \subset \FC(\F_0)$.  Lemma~\ref{lem:contractible} completes the proof.
\endproof

We are now ready for the general case.

\noindent{\bf Proof of  Theorem~\ref{thm:main}}  Let  $m= 3k(2 - 2 \chi(S)) \ge 3k|S_-|$.    Corollary~\ref{reducing sets exist} implies the existence of $K$ and $\cR$ satisfying our standing hypothesis and with the additional feature that  $K$  is $\FC(f_l)$-complete for all $f_l$.  Lemma~\ref{must intersect} implies that for each $l$ there exists $M =N_i$ or $M=A_j$ such that $\FC(f_l) \cap K \cap M \ne \emptyset$.  If there exists a single such $M$ and at least $k$ values of $l$ such that $\FC(f_l) \cap M \cap K \ne \emptyset$ then Lemma~\ref{many fixed points} implies that $\Fix_c(\F) \ne \emptyset$ and we are done.
We may therefore assume that there are at least  $2k|S_-|$ distinct values of $l$ such that $\FC(f_l) \cap K \cap X_M = \emptyset$ for all $M = N_i$ and all  $M=A_j$. After reordering the $f_l$'s we may assume that $\{f_1,\ldots,f_{m'}\}$ have this property where $m' = 2k|S_-|$.

Let  $K'$ be the subset of  $K$ obtained by removing $K \cap X_M$ for all $M$.  Then $K'$ has locally constant Nielsen behavior and is $\FC(f_l)$ complete for $1 \le l \le m'$.   Let $\cR'$ be a TNF-complete normal form for $\{f_l\}$ that contains $\cR$.   (Removing punctures  may create new finitely punctured  components of $S \setminus \cR$  and so may cause $\cR$ to be TNF-incomplete with respect to $K'$.)   Define subsurfaces $N_i'$, $A_j'$ and $S_-'$ using $K'$ and $\cR'$.  
      Note that  $S_-' \subset S_-$ and that each non-disk component of $S_-' \setminus \cR'$ is contained in a non-disk component of $S_- \setminus \cR$.  It follows that  $ K' \cap X_{M'} = \emptyset$ for all $M' = N_i'$ and all  $M'=A_j'$.  

Since $m' = 2k|S_-|$ there exists a single $M'$ and at least $k$ distinct values of $l$ such that 
  $K \cap \FC(f_l) \cap D_l \ne \emptyset$ for some component $D_l$ of $M \setminus X_M$. The fact that $\cR'$ is TNF-complete implies that  each $f_l|X_M$ is isotopic rel $K'$ to a homeomorphism $F_l$ whose restriction to $X_M$ is either the identity, is periodic with period greater than one or is pseudo-Anosov.

If each $F_l|X_M$ is the identity then $K' \cap M \cap \Fix(f_l) \subset \FC(f_l)$ for each $l$ and Proposition~\ref{prop:lift} then completes the proof.    
We may therefore assume that some $F_l$, which we denote simply as $F$,  is either periodic with period greater than one or is pseudo-Anosov.    Every component of $\partial M$ and every component of $\partial X_M$ bounding a disk that intersects $K \cap \Fix(f)$ is $F$-invariant.   Thus $X_M$   has at least three invariant boundary components and $F|X_M$ is not  periodic with period greater than one.   

It remains to consider the case that $F|X_M$  is pseudo-Anosov.  As in the  proof of Lemma~\ref{many fixed points}, there is a finite index subgroup $\F_0$ of $\F$  whose image in the mapping class group of $X_M$ rel $K'$ is the cyclic group generated by $[F|X_M]$.  Thus each disk component of $M \setminus X_M$ that has non-trivial intersection with  $K \cap \FC(f)$ satisfies the hypotheses of  Proposition~\ref{prop:lift}, which then completes the proof. 
\endproof


\begin{thebibliography}{100}
\bibitem{B1}
{\bf C. Bonatti,}
\newblock {\em Un point fixe common pour des diffeomorphisms commutants de $S^2$}
\newblock {Ann. Math} {\bf 129} (1989),  61--79

\bibitem{B2}
{\bf C. Bonatti,}
\newblock {\em Difféomorphismes commutants des surfaces et stabilité des fibrations en tores.}
\newblock {Topology} {\bf 29} (1990), no.2  205--209


\bibitem{Brown}
{\bf R. F. Brown}
\newblock {\em The Lefschetz fixed point theorem.}
\newblock{Scott, Foresman and Co., Glenview, Ill.-London} (1971) vi+186 pp.

\bibitem{Dold}
{\bf A. Dold,}
\newblock {\em Fixed point index and fixed point theorem for Euclidean neighborhood retracts.}
\newblock {Topology} {\bf 4} (1965),  1--8

\bibitem{FLP}
{\bf A.~Fathi, F.~Laudenbach and V.~Poenaru,}
\newblock {\em Travaux de Thurston sur les surfaces,}
{Asterisque} (1979) 66--67

\bibitem{F}
{\bf S. Firmo}
\newblock {\em A note on commuting diffeomorphisms of surfaces.}
\newblock {Nonlinearity} {\bf 18} (2005), n0.4  1511--1526

\bibitem{fh:periodic}
{\bf J.~Franks and M.~Handel,}
\newblock {\em  Periodic points of Hamiltonian surface diffeomorphisms,}
\newblock  {Geom. Topol.} {\bf 7} (2003) 713--756


\bibitem{FHP}
{\bf J.~Franks, M.~Handel and K.~Parwani}
\newblock {\em Fixed Points of abelian actions on $S^2$.}
\newblock{preprint}

\bibitem{han:commuting}
{\bf M.~Handel,}
\newblock {\em Commuting homeomorphisms of $S^2$,}
\newblock {Topology} {\bf 31} (1992) 293--303

\bibitem{HT}
{\bf M.~Handel and W.~Thurston}
\newblock {\em New Proofs of Some Results of Nielsen,}
\newblock {Adv. in Math.} {\bf 56} (1985) 173--191

\bibitem{Th}
{\bf W.~Thurston,}
\newblock {\em On the geometry and dynamics of diffeomorphisms of surfaces,}
\newblock Bull. Amer. Math. Soc. {\bf 19} (1988) 417--431


\end{thebibliography}
\end{document}